\documentstyle[12pt,righttag]{amsart}
\numberwithin{equation}{section}
\theoremstyle{plain}
  \newtheorem{theorem}{Theorem}[section]
  
  \newtheorem{proposition}[theorem]{Proposition}
  \newtheorem{lemma}[theorem]{Lemma}
  \newtheorem {corollary}[theorem]{Corollary}
\theoremstyle{definition}
  \newtheorem{definition}[theorem]{Definition}
\theoremstyle{remark}
 \newtheorem{remark}[theorem]{Remark}
  \newtheorem{example}[theorem]{Example}

\def\e{{\epsilon}}
\def\ech{\epsilon^\vee}
\def\gl{{\mathfrak{gl}}}
\def\h{{\mathfrak h}}

\def\t{{{\mathfrak h}}}

\def\lm{{\lambda}}

\def\gl{{\mathfrak gl}}
\def\al{\alpha}
\def\Sym{{\mathfrak S}}
\def\C{{\Bbb C}}

\def\Z{{\Bbb Z}}

\def\C{{\Bbb C}}
\def\In{{\mathcal I}}
\def\+{\mathop{\oplus}}
\def\*{\mathop{\otimes}} 
\def\aff{\dot}
\def\daff{\ddot}

\def\Rep{{{\mathcal O}}}
\def\pair{{{\lm,\mu}}}
\def\one{{{\bf 1}}}
\def\ch{^\vee}

\def\bra{{\langle}}
\def\ket{{\rangle}}
\def\qed{\hfill$\square$}

\def\dim{{\hbox{\rm{dim}}}\,}

\def\gen{{{\rm{gen}}}}

\def\Irr{{\rm{Irr}}}

\def\aa{{{a}}}
\def\bb{{{b}}}
\def\ir{{K}}
\def\affh{{{\aff\h}}}
\def\gm{{{y}}}
\def\xxi{{{\xi}}}
\def\varpip{{{\varpi_p}}}

\def\Set{{\mathcal S}}
\def\affh{{{\tilde\h}}}
\def\udl{\underline}

\def\tp{{{t}}}
\def\dom{{{\mathcal D}}} 
\def\FDaff{{{\aff\In_{p,\kappa}^{*+}}}}
\def\zFDaff{{{\aff\In_{p,\kappa}^{+}}}}
\def\HFDaff{{{\aff\dom_{p,\kappa}}}}

\def\BB{{{\mathcal B}}}
\def\CS{{{\mathcal C}}}
\newcommand{\brac}[2]{\langle #1 \mid #2 \rangle}
\begin{document}
\title[Degenerate Double Affine Hecke Algebras]{Classification 
of simple modules over degenerate double Affine Hecke algebras 
of type $A$}
\author{Takeshi Suzuki}
\address{Research Institute for Mathematical Sciences\\ 
Kyoto University}
\email{takeshi@@kurims.kyoto-u.ac.jp}
\maketitle
\begin{abstract}
We study a class of representations over
the degenerate double affine Hecke algebra of $\gl_n$
by an algebraic method.
As fundamental objects in this class, 
we introduce certain induced modules and 
study some of their properties.
In particular, it is shown that these induced modules
 have  unique simple quotient modules
under certain conditions.
Moreover, 
we show that any simple module in this class is obtained 
as such a simple quotient,
and give a classification of all the simple modules.
\end{abstract}
\begin{center}
  {\sc Introduction}
\end{center}
Double affine Hecke algebras and their degenerate 
(or graded) version
are introduced by Cherednik~\cite{Ch;ellip},
and successfully applied to 
the theory of symmetric 
polynomials~\cite{Ch;Macdonald,Ch;Macdonald2}.

The purpose of this paper is to give an algebraic
approach to the study of the representation theory 
of  the degenerate
double affine Hecke algebra of type $A$.
In particular, we give a classification of simple modules
of a certain class, which is studied in~\cite{BEG} 
for double affine Hecke algebras
from the geometric viewpoint.

Let $\daff H_n$ denote the degenerate affine Hecke
algebra of $\gl_n$.
The algebra $\daff H_n$ has a commutative subalgebra 
$S(\aff\t)$.
Here $S(\aff\t)$ denotes the symmetric algebra of the 
vector space $\aff\h=\bigoplus_{i=1}^n\C\ech_i\+\C c$ with
$c$ central in $\daff H_n$.
Note that a locally $S(\aff\t)$-finite  $\daff H_n$-module
 admits generalized weight space decomposition with respect to 
the action of $\aff\t$.
We study the category $\Rep(\daff H_n)$ 
consisting of  
finitely generated, locally $S(\aff\t)$-finite $\daff H_n$-modules 
whose (generalized) weights 
are integral.

Since $c$ is a center, it follows that
the category is decomposed into a direct sum of subcategories
$\Rep_\kappa(\daff H_n)\ (\kappa\in\Z)$,
where $\Rep_\kappa(\daff H_n)$ denotes the full subcategory
consisting of modules on which $c$ acts as 
a scalar multiple by $\kappa$.

We will give a classification of all simple modules in
$\Rep_\kappa(\daff H_n)$ with $\kappa\neq 0$.
(We do not treat the case $\kappa=0$, which is rather special.)
Let us sketch our approach.

We introduce a certain  set of parameters $\Set$.
For each parameter $(\pair)\in \Set$, we introduce
an $\daff H_n$-module $\daff M(\pair)$, which is
 induced from a certain one-dimensional module
of a parabolic subalgebra of $\daff H_n$, and 
investigate their properties (\S~\ref{sec;induced}).

The first main result in this paper is
Theorem~\ref{th;usqd},
which gives a sufficient condition on $(\pair)$
ensuring that $\daff M(\pair)$ 
has a unique simple quotient module, denoted by
 $\daff L(\pair)$.
Define $\Set^+$ as the subset of $\Set$
consisting of all parameters satisfying
the conditions in Theorem~\ref{th;usqd}.
Then, we can get a correspondence
from $\Set^+$ to 
the set of isomorphism classes of
simple objects in $\Rep_\kappa(\daff H_n)$.

We prove that
any simple module in $\Rep_\kappa(\daff H_n)$ can be
obtained as a simple quotient $\daff L(\pair)$ 
for some 
$(\pair)\in \Set^+$ (Theorem~\ref{th;allirrd}), that is,
 the correspondence above is surjective.
Furthermore, we write down when two parameters in $\Set^+$
give isomorphic simple modules (Theorem~\ref{th;isomd}).
This completes
the classification of simple objects in $\Rep_\kappa(\daff H_n)$.
It turns out that the  set
 of isomorphism classes of
simple objects in $\Rep_\kappa(\daff H_n)$
is indexed by isomorphism classes of $n$-dimensional nilpotent 
representations of the (cyclic) quiver 
of type $A_{\kappa}^{(1)}$.

We give detailed proofs for all the statements above
by an algebraic and rather direct method 
with the help of some fundamental results
on the representation theory of
the (degenerate) affine Hecke algebra.

We treat the degenerate double affine Hecke algebra 
in this paper, but 
it is easy to modify the arguments to obtain 
the same results
for the double affine Hecke algebra of $\gl_n$
provided that a certain parameter (often denoted by $q$)
is not a root of one. 
Note that 
the classification of simple modules over 
the affine Hecke algebra
can be deduced to the same problem for 
the degenerate affine Hecke 
algebra and vise versa~\cite{Lu},
but 
the corresponding rigorous statement has not been established
(as far as the author knows).

A similar class of representations over 
(non-degenerate) double affine Hecke algebras of general type
 is studied by a geometric method in the preprint \cite{Va},
where the classification of simple modules
 and certain Jordan-Holder multiplicity
formulas are obtained by means of
the theory of perverse sheaves and equivariant K-theory.
In particular, Vasserot's result gives a geometric proof of
our classification for the double affine Hecke algebra.
It should be also mentioned that Cherednik 
announces the classification of simple modules
over the double affine Hecke algebra of type $A$
by an alternative algebraic approach
in the preprint~\cite{Ch;mathQA}.

\medskip\noindent
{\bf Acknowledgment.}
The author would like to thank T.~Arakawa and A.~Tsuchiya for
the collaboration in \cite{AST} and \cite{AS},
which lead to this work.
\section{Affine root system }
Through this paper, we use the notation
$$[i,j]=\{i,i+1,\dots,j\}$$
for $i,j\in\Z$ with $i<j$.

Fix $n\in \Z_{>0}$.
Let 
$\affh$
be an $(n+2)$-dimensional vector space over $\C$
with the basis $\{\ech_1,\ech_2,\dots,\ech_n,c,d\}$:
$\affh=\bigoplus_{i=1}^n \C\ech_i \+\C c\+\C d$.

Introduce the non-degenerate  symmetric bilinear form $(\ |\  )$ 
on $\affh$ by
\begin{align*} 
(\ech_i|\ech_j)&=\delta_{ij},\quad 
(\e\ch_i| c)=(\e\ch_i|d)=0, \\
(c|d)&=1,\quad (c|c)=(d|d)=0.
\end{align*}
Put
$\h=\bigoplus_{i=1}^n \C\ech_i$
and ${\aff\t}=\h\+\C c$.

Let $\affh^*=\bigoplus_{i=1}^n\e_i\oplus \C c^* \oplus \C \delta$ 
be the dual space of $\affh$,
where $\e_i$,  $c^* $ and $\delta$
are   the dual vectors of $\ech_i$,
$c$ and $d$ respectively.

We  identify
the dual space $\aff\t^{*}$ 
(resp. $\h^*$) of $\aff\t$ (resp. $\h$) as a subspace 
of $\affh^{*}$ via the identification
${\aff\t}^{*}={\affh}^{*}/\C \delta\cong \h^{*}\oplus
\C c^{*}$
(resp.
$\h^*=\affh^*/(\C c^*\+\C\delta)\cong
\bigoplus_{i=1}^n\C \e_i$).

The natural pairing is denoted by
$\bra\, |\,\ket : \affh^*\times \affh\to\C$. 
There exists an isomorphism
between $\affh^{*}$ and
$\affh$
such that 
$ \e_i\mapsto\ech_i$,
$\delta\mapsto c$ and $c^{*} \mapsto d$.
We denote by $\zeta\ch\in\affh$ the image of
$\zeta\in\affh^{*}$ under this isomorphism.

Put $\al_{ij}=\e_i-\e_j\ (1\leq i\neq j\leq n)$ and 
$\al_i=\al_{ii+1}\ (1\leq i\leq n-1) $. 
Then 
\begin{align*}
R&=\left\{ \alpha_{ij}\mid
1\leq i\ne j \leq n \right\}
,\
R^+=\left\{ \alpha_{ij}\mid 1\leq i<j \leq n\right\},\\ 
\Pi&=\left\{\alpha_1,\dots,\alpha_{n-1}\right\}
\end{align*}
give the system of roots,
positive roots and simple roots of type $A_{n-1}$
respectively.

Put $\al_0=-\al_{1n}+\delta$, and 
define the system $\aff R$ of roots,  $\aff R^+$
of positive roots and  $\aff \Pi$ of simple roots 
of type $A_{n-1}^{(1)}$ by
$$
\begin{array}{l}
\aff R=
\left\{\alpha+k\delta\,|\,
\alpha\in R,\,k\in\Z\right\},\\
\aff R^+=\left\{\alpha+k\delta~|~
\alpha\in R^+,\,k\in \Z_{\geq0}\right\}
\sqcup \left\{-\alpha+k\delta~|~
\alpha\in R^+,\,k\in\Z_{>0}\right\},\\
\aff\Pi=
\left\{\al_0,\al_1,\dots,\al_{n-1}
\right\}.
\end{array}
$$
\section{Affine Weyl group}
Let $Q$ denote the root lattice 
$\bigoplus_{i=1}^{n-1}\Z \al_i$ and 
let $P$ denote the weight lattice 
$\bigoplus_{i=1}^n \Z\e_i$ of $\gl_n$.
Let $W_n$ denote the Weyl group of $\gl_n$,
which is isomorphic to the symmetric group
${\mathfrak S}_n$.

{The extended affine Weyl group} $\aff W_n $
(resp. {the affine Weyl group} $\aff W_n^\circ$)
of $\widehat\gl_n$ is defined as the semidirect product
of $W_n$ and $P$ (resp. $Q$)
with the relation $w\cdot t_{ \eta}\cdot {w} ^{-1}=
t_{{w}(\eta)}$,
where $w$ and $t_{ \eta} $
are the elements in $\aff W$ corresponding to $w\in W$ and 
$\eta\in P$ (resp. $Q$).
In the following, we simply denote
$\aff W=\aff W_n,\ \aff W^\circ= \aff W^\circ_n$ and 
$W=W_n$.

Let $s_{\al}\in W$ denote the reflection
corresponding to $\al\in R $.
For an affine root
$\beta=\al+k\delta\in \aff R\ (\al\in R ,\ k\in\Z)$, 
define the corresponding affine reflection by
$s_{\beta}=t_{-k\al}\cdot s_{\al}$. 

Put
$s_i=s_{\alpha_i}$
for $i\in [0,n-1]$ and
put $\pi=t_{\e_1}\cdot s_1\cdots s_{n-1}$. 
The following fact is well-known.
\begin{proposition}
${\rm{(i)}}$ The group $\aff W$ is isomorphic to 
the group defined by the following
generators and relations$:$
$$
\begin{array}{ccl}
\hbox{generators}&:& s_i\ (i\in[0,n-1]\cong \Z/n\Z),
\ \pi^{\pm 1}.\\
\hbox{relations}
&:& s_i^2=1,\\
&& s_i s_{i+1}\ s_i=
s_{i+1} s_i s_{i+1} \ (i\in\Z/n\Z),\\
&& s_i s_{j}=s_{j} s_i\ 
(i-j\not\equiv\pm 1 \bmod n),\\
&& \pi s_i=s_{i+1} \pi\ (i\in\Z/n\Z),\\
&& \pi\pi^{-1}=\pi^{-1}\pi=1.
\end{array}
$$
${\rm{(ii)}}$ The subgroup $\aff W^\circ$ is generated by
the simple reflections $s_0,s_1,\dots,s_{n-1}$. 
\end{proposition}

The action of  $\aff W $
on 
$\affh$ 
is given by
the following formulas:
\begin{equation}\label{eq;ac}
\begin{array}{l}
s_{\al} \left( h\right)=
 h-\bra\al| h\ket \al\ch\ \
(\al\in \aff R,\ h\in\affh),\\
\pi(\ech_i)=\ech_{i+1}\ \ (i\in[1,n-1]),
\quad \pi(\ech_n)=\ech_1-c,\\
 \pi(c)=c,\ \ \pi(d)=d.
\end{array}
\end{equation}
It follows that the action of $t_\eta$ $(\eta\in P)$ 
is given by
\begin{align*}
&t_{ \eta } ( h)= h+
\bra\delta| h\ket \eta\ch-
\left( \bra \eta| h\ket+
{1\over 2}( \eta| \eta)\bra\delta| h\ket\right)c.
\end{align*}

The dual action on $\affh^*$ is given by 
\begin{equation*}\label{eq;acdual}
\begin{array}{l}
s_{\al}(\zeta)=
\zeta-(\al|\zeta) \alpha\ \ 
(\al\in \aff R,\ \zeta\in\affh^*),\\
t_{ \eta } (\zeta)=\zeta+
(\delta|\zeta) \eta-
\left( (\eta|\zeta)+{1\over 2}( \eta| \eta)(\delta| \zeta)
\right)\delta\ \ (\eta\in P,\  \zeta\in\affh^* ),\\
\pi(\e_i)=\e_{i+1}\ (i\in[1 ,n-1]),
\quad \pi(\e_n)=\ech_1-\delta,\\
\pi(c^*)=c^*,\ \ \pi(\delta)=\delta.
\end{array}
\end{equation*}
With respect to these actions, the inner products 
on $\affh$ and $\affh^*$ are $\aff W$-invariant.

Note that the subspace  $\aff\t=\h\+\C c$ is preserved by $\aff W$,
and that the dual action of $\aff W$
on ${\aff\t}^*$ (called the affine action) is given by
\begin{equation}\label{eq;affineaction}
\begin{array}{l}
s_{\al} \left(\zeta\right)=
\zeta-( \al| \zeta) \al,\\
t_{ \eta} (\zeta)=\zeta+(\delta|\zeta)  \eta
\end{array}
\end{equation}
for  $\zeta\in\aff\t^*$, $\al\in\aff R$ and $\eta\in P$.

For $w\in\aff W $,
set 
$$ R(w)=\aff R^+\cap w^{-1}\aff R^-,$$
where $\aff R^-=\aff R\setminus \aff R^+$.
The length $l(w)$ of $w\in\aff W$ is defined as the number  
$\sharp R(w)$
of the elements in $ R(w)$.
For $w\in\aff W $,
an expression $w=\pi^k\cdot s_{j_1}\cdots s_{j_m}$ 
is called a reduced expression
if $m=l(w)$.
It can be seen that
\begin{equation}
  R(w)=\{s_{j_m}\cdots s_{j_2}(\al_{j_1}),
s_{j_m}\cdots s_{j_3}(\al_{j_2}),\dots,\al_{j_m}\}
\end{equation}
if  $w=\pi^k\cdot s_{j_1}\cdots s_{j_m}$ is a reduced 
expression.

The partial ordering $\preceq $ is 
defined
in the Coxeter group $\aff W^\circ$ as follows:
$w\preceq w'$
$\Leftrightarrow$ $w$ can be obtained as a subexpression
of a reduced expression  of $w'$. 
Extend this ordering $\preceq $ to
the  partial ordering  in $\aff W$
by 
$\pi^k w\preceq \pi^{k'}w'\Leftrightarrow k= k' 
\hbox{ and } 
w\preceq w'
$
($k,k'\in\Z,\ w,w'\in \aff W^\circ$).

\medskip
Let $I$ be a subset of  $[0,n-1]$.
Put 
\begin{align*}
\aff\Pi_I&=\{\al_i\mid i\in I\}\subseteq\aff\Pi,\\
\aff W_I&=\bra s_i~;~ i\in I\ket\subseteq \aff W,\\
\aff R_I&=\{\al\in\aff R\mid s_\al\in \aff W_I\}.
\end{align*}
Note that $\aff W_I$ is the parabolic subgroup corresponding 
to $\aff\Pi_I$.
Define 
\begin{equation*}
 \aff W^{I}=\left\{w\in \aff W\mid R(w)\subset 
\aff R^+\setminus(\aff R^+\cap \aff R_I)\right\}.
\end{equation*}
The following fact is well-known.
\begin{proposition}
\label{pr;coset}
For any $w\in\aff W$, 
there exist a unique $w_1\in \aff W^{I}$ 
and  a unique $u\in \aff W_{I}$, such that $w=w_1\cdot u$. 
Their length satisfy $l(w)=l(w_1)+l(u)$.
In particular, the set $\aff W^{I}$
gives a complete set of  
representatives in the coset $\aff W/ \aff W_{I}$.
\end{proposition}
In the case $I\subseteq [1,n-1]$, 
we can define $W_I(=\aff W_I)$ and $W^I\subseteq W$
analogously,
and similar statements as Proposition~\ref{pr;coset} 
hold for them.

Put
$P^-=\{\xi\in P\mid (\xi\mid \al)\leq 0
\hbox{ for any }\al\in R^+\}.$

\begin{lemma}\label{lem;gamma}
Let $\eta\in P$. Let $y$ be a shortest element
 of
$W$ such that $y(\eta)\in P^-$.
Then 
$R(y)=\{\al\in R^+\mid (\eta\mid \al)>0\}.$
\end{lemma}
\noindent{\it Proof.}
Let $\al\in R({y})=R^+\cap {y}^{-1}R^-$.
Then $(\eta\mid \al)=(y(\eta)\mid y(\al))
\geq 0$.
If $(\eta\mid \al)=0$ then $s_{\al}(\eta)=\eta$. Hence
we have $y s_\al(\eta)\in P^-$ and $l(ys_\al)<l(y)$.
This contradicts to the choice of $y$.
Therefore we have $(\eta\mid \al)>0$ and hence
 $R({y})\subseteq\{\al\in R^+\mid (\eta\mid \al)>0\}$. 
It is easier to show the opposite inclusion.\qed

\medskip
We denote by  $\gm_{\eta}$ 
the (unique) shortest element of
$W$ such that $\gm_{\eta}(\eta)\in P^-$.

The following proposition follows from 
Proposition~\ref{pr;coset}
and Lemma~\ref{lem;gamma}. We omit the detailed proof.
\begin{lemma}$(\cite{AST})$  
\label{lem;representative}
${\rm(i)}$ We have
$
\aff W^{[1,n-1]}=\left\{ t_\eta\cdot
\gm_{\eta}^{-1}\mid \eta\in P\right\}.
$

\smallskip\noindent
${\rm(ii)}$ 
For a subset $I\subset [1,n-1],$ we have
$$
\aff W^{I}=\aff W^{[1,n-1]}\cdot  W^{I}.$$
Moreover, $l(w)=l(t_{\eta}\cdot \gm_{\eta}^{-1})+l(u)$
for $w=t_{\eta}\cdot \gm_{\eta}^{-1}\cdot u$
$(\eta\in P$, $u\in  W^{I})$.
\end{lemma}
\section{Degenerate double affine Hecke algebra}
\label{ss;DDAHA}
Let $\C[\aff W]$ denote the group algebra of $\aff W$ and
let $S(\aff \t)$ denote the symmetric algebra 
of $\aff \t=\h\+\C c$.

The degenerate double affine Hecke algebra was 
introduced by Cherednik 
\cite{Ch;ellip}. 
\begin{definition}
\label{Hdef}
The degenerate double affine Hecke algebra 
$\aff H_n$ of $\gl_n$ is the unital
associative $\C$-algebra
defined by the following properties:

\smallskip\noindent
\rm{(i)} As a $\C$-vector space,
$$
\daff H_n=\C[\aff W]\otimes S(\aff \t).
$$

\smallskip\noindent
\rm{(ii)} The natural inclusions
$\C[\aff W]\hookrightarrow \daff H_n$
and $S(\aff \t)\hookrightarrow \daff H_n$
are algebra homomorphisms
(the images of $w\in\aff W$ and $ h\in\aff\t$ will be
simply denoted by $w$ and $ h$).

\smallskip\noindent
\rm{(iii)} The following relations hold in $\daff H_n${\rm :}
\begin{eqnarray}
\label{eq;rel}
& &s_\al  h-s_\al
 ( h) s_\al=-\bra\al| h\ket\quad 
(\al\in \aff R,\, h\in\aff\t),
\\ 
\label{eq;rel2}
& &\pi  h=\pi( h)\pi\quad ( h\in \aff\t). 
\end{eqnarray}
\end{definition}
By definition, the element $c\in\daff H_n$ belongs to the center 
${\mathcal Z}(\daff H_n)$ of $\daff H_n$.
For $\kappa\in\C^*$, we set
$\daff H_n(\kappa)=\daff H_n/\langle c-\kappa\rangle$.

\begin{definition}
Define the {degenerate affine Hecke algebra} 
${\aff H_n}$
as the subalgebra of $\daff H_n$ generated by
the elements in $W$ and the elements in $\h$:
$$\aff H_n
=\C[W]\otimes S(\h)\subset\daff H_n. $$
\end{definition}
\begin{proposition}
\label{pr;relation}
For $w\in\aff W$ and $ h\in\aff\t$, we have
$$ h w=
w\left(w^{-1}( h)+\sum_{\al\in R(w)}
\bra w(\al)| h\ket s_{\al}\right).
$$
In particular,
$ h  w=w\cdot w^{-1}( h)+
\sum_{w'\prec w}c_{w'}w'$
for some $c_{w'}\in\C$.
\end{proposition}
\noindent
{\it Proof.}
The statement is shown by
by the induction on $l(w)$ using the fact that
 $R(s_iw)=R(w)\sqcup\{w^{-1}(\al_i)\}$
if $l(s_iw)=l(w)+1$.\qed

\medskip
It is easy to verify the following proposition directly
(see e.g. \cite{AST,Lu}).
\begin{proposition}\label{pr;center}
{\rm{(i)}} The center 
of
$\daff H_n$ is given by ${\mathcal Z}(\daff H_n) =\C[c]$.

\smallskip\noindent
{\rm{(ii)}}
The center 
of $\aff H_n$ is given by
$${\mathcal Z}(\aff H_n)=\{\xi\in S(\h)\mid w(\xi)
=\xi\text{ for all }w\in W\}.$$
\end{proposition}
For $i\in\Z$, we introduce
the following notations:
\begin{align}
 \e_i&=\e_{\udl{i}}-k\delta\in \affh^*,\quad
\e\ch_i=\e\ch_{\udl i}-kc \in \affh,
\end{align}
where  $i=\udl{i}+kn$ with
$\udl{i}\in[1,n]$ and $k\in \Z$.

Put 
$\al_{ij}=\e_i-\e_j$ 
(and $\al\ch_{ij}=\e\ch_i-\e\ch_j$) for any $i,j\in\Z$.
Note that
\begin{equation*}
\begin{array}{ll}
  \al_{ij}\in \aff R\ &\Leftrightarrow
\ i\not\equiv j \text{ mod }n,\\
 \al_{ij}\in \aff R^+\ &\Leftrightarrow
\ i\not\equiv j \text{ mod }n \text{ and }i<j.
\end{array}
\end{equation*}

Let $J=\{j_1,j_2,\dots,j_m\}$ be a subset of $\Z$
such that $\al_{j_a j_b}
=\e_{j_a}-\e_{j_{b}}\in \aff R$
for $a\neq b$.
Then, in particular, we have $m\leq n$.

Define $\aff H_J$ to be the subalgebra of $\daff H_n$
generated by
$$
\ech_{j_1},\ech_{j_2}
,\dots,\ech_{j_m},\ 
s_{\al_{j_1 j_2}}, s_{\al_{j_2 j_3}},\dots,s_{\al_{j_{m-1} j_m}}.
$$
The following lemma will be used later.
\begin{lemma}
\label{lem;Jsub}
Let 
$J=\{j_1,j_2,\dots,j_m\}$ be a subset of $\Z$
such that $j_a\not\equiv j_b$ mod $n$ for $a\neq b$.
Suppose that $j_1<j_2<\cdots<j_m$.
Then, the algebra $\aff H_J$ is
isomorphic to the degenerate affine Hecke algebra $\aff H_m$
of $\gl_m$.
\end{lemma}
\noindent{\it Proof.}
Note that
$\al_{j_a j_{a+1}}
\in \aff R^+$ for all $a\in[1,m-1]$.
Using Proposition~\ref{pr;relation},
it is  verified that
there exists an algebra homomorphism
$\aff H_m\to \aff H_J$
such that
$\ech_i\mapsto \ech_{j_i}$ $(i\in[1,m])$
and $s_i\mapsto s_{\al_{j_i j_{i+1}}}$ $(i\in[1,m-1])$.
This gives an isomorphism.\qed
\begin{example}\label{ex;HJ}
We have $\aff H_{[1,n]}=\aff H_n
$ by definition.
More generally, we have
$\aff H_{[j,j']}\cong \aff H_{j'-j+1}$
for any $j,j'\in\Z$ such that $j'-j\in[1,n-1]$.
\end{example}
\section{Set of parameters and affine Weyl group}
We introduce some sets of parameters, which will index
representations of $\daff H_n.$

Fix $p\in\Z_{>0}$.
Put $C_p(n)=\{(r_1,\dots,r_p)\in 
(\Z_{\geq 0})^p\mid \sum_{i=1}^p r_i=n\}$, and
set 
\begin{align}
\label{eq;I_p}
\In_p&=
\{(\lm,\mu)\in \Z^p\times \Z^p\mid \lm-\mu\in C_p(n)\},\\
\label{eq;I_past}
\In_p^*&=
\{(\lm,\mu)\in \In_p\mid \lm_i-\mu_i>0\text{ for all }i\in[1,p]\}.
\end{align}
We denote 
the extended affine Weyl group of $\gl_p$ by 
$\aff\Sym_p$ instead of $\aff W_p$ in order to avoid confusion.
The elements of $\aff\Sym_p$
corresponding to $s_i,\pi$ and $t_{\e_i}$
are denoted by $\sigma_i,\varpip$
and $\tp_{e_i}$ respectively.
($\{e_i\}_{i\in[1,p]}$ is the generators of 
the weight lattice $P_p$ of $\gl_{p}$:
$P_p=\oplus_{i=1}^p\Z e_i$.)
We put $\aff \Sym_1=\bra \varpi_1\ket$ for convenience.
The subgroups corresponding to $\aff W_p^\circ$ and $W_p$ are denoted
by  $\aff\Sym_p^\circ$ and $\Sym_p$ respectively:
$$ 
\aff \Sym_p^\circ=
\bra \sigma_0,\sigma_1,\dots,\sigma_{p-1}\ket,\quad
 \Sym_p=\bra \sigma_1,\dots,\sigma_{p-1}\ket.
$$
For $\kappa\in\C$, there exists  an action of $\aff\Sym_p$ 
on the set $\C^p$ which is given by
\begin{align*}
&\sigma_i\circ \lm=
(\lm_1,\dots,\lm_{i+1}-1,\lm_{i}+1,\dots,\lm_p)
\ \ (i\in[1,p-1]),
\label{eq;actionsigma_i}\\
&\sigma_0\circ
\lm=
(\lm_p+\kappa-p+1,\lm_2,\dots,\lm_{p-1},
\lm_1-\kappa+p-1),\label{eq;actionsigma_0}
\\
&\varpip\circ
\lm=
(\lm_p+\kappa-p+1,\lm_1+1,\dots,
\lm_{p-1}+1)
\label{eq;actionpi}
\end{align*}
for $\lm=(\lm_1,\lm_2,\dots,\lm_p)\in\C^p$.
It follows that the action of $\tp_{e_i}$ is given by
\begin{equation*}
\tp_{e_i}\circ\lm=
(\lm_1,\dots,\lm_{i}+\kappa,\dots,\lm_p)\ \ (i\in[1,p]).
\end{equation*}
If $\kappa\in \Z$, then this action preserves $\Z^p$ and 
induces an action of $\aff\Sym_p$ on $\In_p$ and $\In_p^*$
via 
\begin{equation}\label{eq;dotaction2}
w\circ(\pair)=(w\circ\lm,w\circ\mu).
\end{equation}
In the following, we always assume $\kappa\in \Z$.

For $\lm\in\C^p$, put 
\begin{align*}
  \label{eq;[myu0]}
[\lm]_0&=\kappa-p+1-\lm_1+\lm_p,\\
\label{eq;[myui]}
[\lm]_i&=\lm_i-\lm_{i+1}+1\ \ (i\in [1,p-1]).
\end{align*}
\begin{remark}
It is natural
to describe the action $\circ$ and the numbers $[\lm]_i$ 
 in terms of the root system of type $A_{p-1}^{(1)}$:
Put $\aff\t_p^*=\bigoplus_{i=1}^p\C e_i\+ \C c^*$,
 where
notations are analogous to the $A_{n-1}^{(1)}$ case,
and regard $\lm\in\C^p$
as  an element of $\aff\t_p^*$ by
$\lm=\sum\lm_i e_i+(\kappa-p)c^*$.
Then  we have 
\begin{align*}
&w\circ\lm=w(\lm+\rho)-\rho\quad (w\in\aff\Sym_p),\\  
&[\lm]_i=\bra\lm+\rho\mid\al\ch_i\ket\quad (i\in[0,p-1]),
\end{align*}
where $\rho=\sum_{i=1}^{p}(-i+1) e_i+p c^*$.
\end{remark}
Set
\begin{align*}
\dom_p &=
\{ \lm\in \Z^p \mid [\lm]_i\geq 0
\hbox{ for all } i\in [1,p-1] \},\\
\HFDaff&=\{ \lm\in \Z^p\mid
[\lm]_i\geq 0
\hbox{ for all } i\in [0,p-1] \}.
\end{align*}
The following fact is well-known:
\begin{lemma}\label{lem;fd1}
Let $\kappa\in \Z_{>0}$.

\smallskip\noindent
$\rm{(i)}$
$\dom_p$ is a fundamental domain 
for the action of $\Sym_p$ on $\Z^p$.

\smallskip\noindent
$\rm{(ii)}$
$\HFDaff$ is
a fundamental domain for 
the action of $\aff\Sym_p^\circ$ on $\Z^p$.
\end{lemma}
The proof of the following lemma is similar to the proof of
 Lemma~\ref{lem;gamma}.
\begin{lemma}
\label{lem;fundshortest}
Let $\kappa\in \Z_{>0}$ and $\lm\in\Z^p$.
Let $w$ be the shortest element in $\aff \Sym_p^\circ$  
such that $w\circ\lm\in \HFDaff$.
Then, we have  $[\sigma_{i_{k+1}}\sigma_{i_{k+2}}\cdots 
\sigma_{i_l}\circ\lm]_{i_k}<0$
for each $k=[1,l]$, where
$w=\sigma_{i_1}\sigma_{i_2}\cdots \sigma_{i_l}$ is 
a reduced expression of $w$.
\end{lemma}

\medskip
For $\mu\in\Z^p$, set
\begin{align*}
  \dom_p^{\mu}&=\{\lm\in\Z^p \mid [\lm]_i\geq 0 
\hbox{ for any }i\in
[1,p-1] \hbox{ such that }[\mu]_i=0\},\\
  \HFDaff^{\mu}&=\{\lm\in\Z^p\mid [\lm]_i\geq 0 
\hbox{ for any }i\in
[0,p-1] \hbox{ such that }[\mu]_i=0\}.
\end{align*}
Put
\begin{align}\label{eq;I_pplus}
\In_p^+ &=
\{ (\pair)\in \In_p\mid \mu\in \dom_p,\
\lm\in \dom_p^\mu\},\\
\label{eq;I_pk}
\zFDaff  &=
\{ (\pair)\in \In_p\mid \mu\in \HFDaff,\ 
\lm\in \HFDaff^\mu\},\\
\In_p^{*+}&=\In_p^*\cap\In_p^+,\ \ 
\FDaff=\In_p^*\cap\zFDaff.
\end{align}
It is easy to show the following lemma:
\begin{proposition}\label{pr;fdpair}
Let $\kappa\in\Z_{>0}$. 

\smallskip\noindent
$\rm{(i)}$
$\In_p^+$ $($resp. $\In_p^{*+})$
is a fundamental domain for 
the action \eqref{eq;dotaction2} of 
$\Sym_p$ on $\In_p$ $($resp. $\In_p^*)$.

\smallskip\noindent
$\rm{(ii)}$
$\zFDaff$ $($resp. $\FDaff)$
is a fundamental domain 
for the action \eqref{eq;dotaction2} of $\aff\Sym_p^\circ$
on $\In_p$ $($resp. $\In_p^*)$.

\smallskip\noindent
$\rm{(iii)}$
$\varpip$ preserves the sets $\zFDaff$ and 
$\FDaff$
respectively.
\end{proposition}
\section{Representations of degenerate affine Hecke algebras}
We review some facts on the representation theory 
of the degenerate affine Hecke algebra $\aff H_n$
for later use.

For an $\aff H_n$-module
 $N$ and $\zeta\in\h^{*}$, define 
the weight space $N_{\zeta}$ 
and the generalized weight space $N_{\zeta}^{\gen}$
of weight $\zeta$ by
\begin{align*}
N_{\zeta}&=\left\{v\in N\mid h v=\bra\zeta| h\ket v 
\hbox{ for any }  h\in \h\right\},\\
N_{\zeta}^{\gen}&=
\left\{v\in N \mid(  h-\bra \zeta| h\ket)^k v=0
\hbox{ for any }  h\in \h, 
\text{ for some }k \in \Z_{>0}\right\}.
\end{align*}
Denote by $P(N)$ the set of all weights of $N$:
$$P(N)=\{\zeta\in\h^*\mid N_\zeta\neq \{0\}\}
=\{\zeta\in\h^*\mid N^\gen_\zeta\neq \{0\}\}.
$$
Define $\Rep(\aff H_n)$ to be the category consisting of all 
finite-dimensional $\aff H_n$-modules $N$ such that
$N=\bigoplus_{\zeta\in P} N_\zeta^\gen$,
i.e. $P(N)\subseteq P$.

Let $(\pair)\in \In_p$.
Put
\begin{align}\label{eq;n_i1}
n_0=0,\quad n_i&=\sum_{j=1}^{i}(\lm_j-\mu_j)\ \ (i=[1,p]).
\end{align}
Set
$I_{\pair}=[1,n-1]\setminus \{n_1,n_2,\dots,n_{p-1}\}$.
Then $W_{I_\pair}
=W_{\lm_1-\mu_1}\times W_{\lm_2-\mu_2}\times
\cdots \times W_{\lm_p-\mu_p}$.
We denote $W_{\lm-\mu}=W_{I_\pair}$ 
and $W^{\lm-\mu}= W^{I_\pair}$.
Define $\aff H_{\lm-\mu}$ 
as the subalgebra of $\daff H_n$ generated by 
the elements in $W_{\lm-\mu}$ and the elements in 
$S(\h)$:
\begin{equation*}
\begin{array}{ll}
\aff H_{\lm-\mu}
&=\C[ W_{\lm-\mu}]\otimes S(\h)\\
&=\aff H_{[n_0+1,n_1]}\otimes \aff H_{[n_1+1,n_2]}\otimes
\cdots\otimes\aff H_{[n_{p-1}+1,n_p]}
\subseteq \aff H_n,
\end{array}
\end{equation*}
where $\aff H_{[i,j]}$ is as in Lemma~\ref{lem;Jsub}.

Define $\zeta_{\lm,\mu}$ to be the element of $\h^*$
such that
\begin{equation}\label{eq;weight_of_one}
\bra\zeta_\pair|\ech_j\ket
=
\mu_i-i+j-n_{i-1}\
\ \text{for}\ j\in [n_{i-1}+1,n_{i}].
\end{equation}
Note, in particular, that we have
\begin{equation}\label{eq;wtedge}
\begin{array}{c}
\bra\zeta_\pair|\ech_{n_{i-1}+1}\ket=\mu_i-i+1,\\
\bra\zeta_\pair|\ech_{n_{i}}\ket=\lm_i-i 
\end{array}
\end{equation}
if $n_{i-1}<n_{i}$.
There exists a one-dimensional representation 
$\C \one_\pair$ 
of 
$\aff H_{\lm-\mu}$ such that
\begin{equation}\label{eq;one_pair}
\begin{array}{ll}
{w}\one_\pair=\one_\pair
\quad &\text{for all }w\in  W_{\lm-\mu},\\  
h
{\one}_\pair=
\bra\zeta_\pair| h\ket{\one}_\pair
\quad &\text{for all } h\in \h.
\end{array}
\end{equation}
Define the induced representation $\aff M(\pair)$ 
of $\aff H_n$ associated with $(\pair)$
by
$$\aff M(\pair)=\aff H_n\otimes_{\aff H_{\lm-\mu}} 
\C \one_\pair.
$$
Clearly,
$
\aff  M(\pair)\cong \C[ W/ W_{\lm-\mu}]$
 as a $W$-module.

The induced module $\aff M(\pair)$ is not irreducible in general.
We will use the following criterion
for the irreducibility in  the case $p=2$.
See e.g. \cite{Zel;induced} or \cite{S,S2}
for the proof.
\begin{lemma}
\label{lem;cpser}
Let $(\lm,\mu)\in \In_2^*$ with $\lm=(\lm_1,\lm_2),\
\mu=(\mu_1,\mu_2)$.

\smallskip\noindent
${\rm (i)}$  $\aff M(\lm,\mu)$ is reducible if and only if
one of the following conditions hold$:$

${\rm (a)}$ $\mu_2\leq \mu_1\leq \lm_2+1$ and $\lm_2\leq\lm_1$.

${\rm (b)}$ $\mu_1 \leq\mu_2\leq \lm_1+1$ and $\lm_1\leq\lm_2$.  

\noindent
In each case, there exist the following exact sequences$:$

\smallskip
${\rm (a)}$ 
$0\to \aff L(\sigma\circ\lm,\mu)\to \aff M(\pair)\to \aff 
L(\pair)\to 0.$

${\rm (b)}$ $0\to \aff L(\sigma\circ\lm,\sigma\circ\mu)\to 
\aff M(\pair)\to \aff 
L(\lm,\sigma\circ\mu)\to 0.$

\smallskip
\noindent
Here $\sigma=\sigma_1\in \Sym_2$.

\smallskip\noindent
${\rm (ii)}$
If $\aff M(\lm,\mu)$ is irreducible, then 
$\aff M(\lm,\mu)\cong \aff M(\sigma\circ\lm,\sigma\circ\mu)$.
\end{lemma}
The following lemma follows from Proposition~\ref{pr;relation}.
\begin{lemma}
\label{lem;affine_decom_gen}
We have
\begin{align*}
  P(\aff M(\pair))&=W^{\lm-\mu}\zeta_\pair:=
\left\{w(\zeta_\pair)\mid w\in W^{\lm-\mu}\right\},\\
{\rm dim} \aff M(\pair)_{\xi}^{{\rm \gen}}&=
\sharp\left\{w\in W^{\lm-\mu}\mid w(\zeta_\pair)=\xi\right\}
\quad \text{for }\xi\in \h^*.
\end{align*}
In particular, we have
${\rm dim}\aff M(\pair)_{\zeta_\pair}^{{\rm \gen}}=
\sharp \left( 
W^{\lm-\mu}\cap W[\zeta_\pair]\right)$,
where 
$W[\xi]=\{ w\in W\mid w(\xi)=\xi\}$ for $\xi\in\h^*.$
\end{lemma}
Let $(\pair)\in \In_p$.
Take integers $a_1<a_2<\cdots<a_k$
such that
\begin{align*}
&\{a_1,a_2,\dots,a_k\}
=\{a\in[1,p]\mid[\lm]_a\neq0\text{ or } [\mu]_a\neq 0\}.
\end{align*}
Put $X_\pair=[1,n]\setminus \{n_{a_1},n_{a_2},\dots,n_{a_k}\}$,
where $n_i$ is as in \eqref{eq;n_i1}.
\begin{lemma}
  \label{lem;block}
Let $(\pair)\in\In_p^{+}$.
Then  
$W^{\lm-\mu}\cap W[\zeta_\pair]\subseteq W_{X_\pair}.
$
\end{lemma}
\noindent{\it Proof.}
Take integers $b_1<b_2<\cdots<b_l$
such that
\begin{align*}
&\{b_1,b_2,\dots,b_l\}
=\{a\in[1,p]\mid[\mu]_a\neq 0\},
\end{align*}
and put $Y_\mu=[1,n]\setminus 
\{n_{b_1},n_{b_2},\dots,n_{b_l}\}$.
Then $W_{\lm-\mu}\subseteq W_{X_\pair}\subseteq W_{Y_\mu}$.

Let $w\in W^{\lm-\mu}\cap W[\zeta_\pair]$.
First, we prove 
$w\in W_{Y_\mu}$.

It is enough to show
$w([1,n_{b_i}])=[1,n_{b_i}]$ for all $i\in[1,l]$.
Suppose $w([1,n_{b_i}])\neq [1,n_{b_i}]$.
Let $m$ be the smallest number such that
$m\in [1,n_{b_i}]$ and $w^{-1}(m)\notin [1,n_{b_i}]$.
Then it follows from $w\in W^{\lm-\mu}$
that $w^{-1}(m)=n_{b-1}+1$ for some $b>b_i$.
Since  $w\in  W[\zeta_\pair]$, we have
$$
\mu_b-b+1=\brac{\zeta_\pair}{w^{-1}(\ech_m)}
=\brac{\zeta_\pair}{\ech_m}\geq
\mu_{b_i}-b_i+1.$$
This implies $[\mu]_{b_i}=0$ as $(\pair)\in\In^+_p$, 
and this is a contradiction.
Therefore we have $w\in W_{Y_\mu}.$

Next, let us see
$w\in W_{X_\pair}$. Suppose 
that $w([1,n_{a_i}])\neq [1,n_{a_i}]$ for some $i\in[1,k]$.
Let $m'$ be the largest number such that
$m'\in [1,n_{a_i}]$ and $w(m')\notin [1,n_{a_i}]$.
Then we have 
 $m'=n_{a}$ for some $a>a_{i}$.
By similar arguments as above, we have
 $[\lm]_{a_i}=0$, and this is a contradiction.
Hence we have $w\in W_{X_\pair}.$
\qed

\medskip\noindent
In~\cite{S2}, Lemma~\ref{lem;block}
is used to reduce the proof of the following 
proposition to the special case
$$\lm=(m,m+1,\dots,m+p-1),\ \mu=(0,1,\dots,p-1)$$
with $mp=n$:
\begin{proposition}$(\rm{Lemma~5.2\ in\ }\cite{S2})$
\label{pr;affine}
Let $(\pair)\in \In_p^+$. Then we have
$\aff M(\lm,\mu)_{\zeta_{\lm,\mu}}
=\C\one_\pair$.
\end{proposition}
As a direct consequence of Proposition~\ref{pr;affine},
we have the following.
\begin{theorem}
$(\cite{Ro,S2})$
\label{th;affine1}
Let $(\lm,\mu)\in \In^+_p$. 
Then 
$\aff M(\lm,\mu)$ has a unique simple quotient,
which we denote by $\aff L(\lm,\mu)$.
\end{theorem} 
Let $\Irr (\Rep (\aff H_n))$ denote the set of isomorphism 
classes of irreducible modules in $\Rep (\aff H_n)$. 
Then, by Theorem~\ref{th;affine1}, we
have correspondences
$\In_p^{*+}\to \Irr (\Rep (\aff H_n))$ $(p\in[1,n])$
given by $(\pair)\mapsto \aff L(\lm,\mu)$. 

The classification of simple modules described below
is originally obtained
by Zelevinsky~\cite{Zel;induced} (see also \cite{Ro})
for the affine Hecke algebra.
An alternative algebraic proof of the classification
using Theorem~\ref{th;affine1} is given in \cite{S2}.
%
\begin{theorem} $(See\ \rm{\cite{S2}}$-$\rm{\S 6.})$
\label{th;paraffine}
The correspondence 
$$\bigsqcup_{p=1}^n \In_p^{*+}\to \Irr (\Rep (\aff H_n))$$
which maps $(\pair)\in\ \In_p^{*+}$ $(p\in[1,n])$
to $\aff L(\pair)$
is a bijection.
\end{theorem}
Recall that $\Sym_p$ acts on $\In_p^*$ and
that $\In_p^{*+}$ is a fundamental domain for this action
(Proposition~\ref{pr;fdpair}):
$\In_p^{*+}\cong \In_p^*/ \Sym_p$.
Hence Theorem~\ref{th;paraffine} asserts that
there exists a one to one correspondence
$$
\bigsqcup_{p=1}^n \In_p^{*}/ \Sym_p
\leftrightarrow\Irr (\Rep (\aff H_n)).$$
\section{Induced representations of $\daff H_n$}
\label{sec;induced}
Let $\kappa\in \Z$.
We consider representations of $\daff H_n(\kappa)$,
namely, representations of 
$\daff H_n$ on which $c\in {\mathcal Z}(\daff H_n)$
acts as a constant integer $\kappa$. 
For an $\daff H_n$-module
 $N$ and a weight $\zeta\in\aff\t^{*}$, 
we use the same notations as those for $\aff H_n$-modules
to denote
the weight space,
the generalized weight space 
and  the set of all weights of $N$:
\begin{align*}
N_{\zeta}&=\left\{v\in N\mid h v=\bra\zeta| h\ket v 
\hbox{ for any }
  h\in \aff\t\right\},\\
N_{\zeta}^{\gen}&=
\left\{v\in N \mid(  h-\bra \zeta| h\ket)^k v=0\hbox{ for any }
  h\in \aff\t, 
\text{ for some }k \in \Z_{>0}\right\},\\
P(N)&=\{\zeta\in\aff\h^*\mid N_\zeta\neq \{0\}\}.
\end{align*}
If $N$ is an $\daff H_n(\kappa)$-module,
then 
any weight of $N$ is of the form
$\zeta+\kappa c^*$ for some $\zeta\in\h^*$.
Put $P_\kappa=P+\kappa c^*\subseteq \aff \t^*$. 
Note that $P_\kappa$ is preserved under the action 
\eqref{eq;affineaction} of $\aff W$.

\begin{definition}
Define $\Rep_\kappa(\daff H_n)$ to be 
the full subcategory of the category
of finitely generated $\daff H_n(\kappa)$-modules
consisting of  
those $\daff H_n(\kappa)$-modules 
$N$ such that 

\smallskip\noindent
(i) $N$ is locally $S(\h)$-finite,

\smallskip\noindent
(ii)  $P(N)\subseteq P_\kappa$.

\smallskip\noindent
(Note the the condition (i) ensures 
the generalized weight space decomposition
$N=\bigoplus_{\zeta\in \aff\t^*}N_{\zeta}^\gen$.)
\end{definition}
\begin{remark}
A similar category 
for the (non-degenerate)
double affine Hecke algebra of general type
is studied from a geometric viewpoint
in \cite{BEG},
where Deligne-Langlands-Lusztig type conjecture
concerning the classification of simple modules
is proposed. 
\end{remark}
\begin{remark}
There exists an algebra automorphism on $\daff H_n$
such that
$$s_i\mapsto -s_i,\ (i\in [0,n-1]),\ \pi\mapsto\pi,\ 
\ \e\ch_i\mapsto -\e\ch_i\ (i\in [1,n]),
\ c\mapsto -c.$$
This gives a categorical equivalence
$\Rep_\kappa(\daff H_n)\cong \Rep_{-\kappa}(\daff H_n)$.
\end{remark}
In the rest of this paper, 
we mostly consider the case $\kappa\in \Z_{>0}.$

For $(\pair)\in \In_p$, 
put $\zeta^\kappa_\pair=\zeta_\pair+\kappa c^*\in
P_\kappa$, where $\zeta_\pair$ is given in 
\eqref{eq;weight_of_one}.
We regard the one-dimensional $\aff H_{\lm-\mu}$-module 
$\C\one_\pair$ (defined by \eqref{eq;one_pair})
as an $(\aff H_{\lm-\mu}\otimes\C [c])$-module 
by letting $c$ act as a constant integer $\kappa$; we have
\begin{equation*}
\begin{array}{ll}
{w} {\one_\pair}={\one_\pair} 
\quad &\text{for all }w\in  W_{\lm-\mu},\\  
h{\one_\pair}=
\bra\zeta^\kappa_\pair| h\ket{\one_\pair}
\quad &\text{for all } h\in \aff\t.
\end{array}
\end{equation*}

Define an $\daff H_n$-module
$\daff M(\pair)$ by
$$\daff M(\pair)=\daff H_n\otimes_{\aff H_{\lm-\mu}\otimes\C [c]}
\C {\one_\pair}.$$
Clearly,
\begin{align*}
  \daff M(\pair)&\cong \daff H_n(\kappa)\otimes_{\aff H_n}
\aff M(\pair)\ \ \hbox{as an }
\daff H_n(\kappa)\hbox{-module},\\
&\cong \C[P]\otimes \aff M(\pair)\ \ \hbox{as an }
\aff H_n\hbox{-module},\\
&\cong \C[\aff W/ W_{\lm-\mu}]\ \ \hbox{as a }
\aff W\hbox{-module}.
\end{align*}
By Proposition~\ref{pr;relation}, 
we have the following:
\begin{proposition}
\label{pr;decom_gen}
We have

\smallskip\noindent
$\rm{(i)}$ 
 $P(\daff M(\pair))=\aff W^{\lm-\mu}\zeta_\pair^\kappa:=
\left\{w(\zeta_\pair^\kappa)\mid w\in 
\aff W^{\lm-\mu}\right\}.$

\smallskip\noindent
${\rm(ii)}$ 
${\rm dim} \daff M(\pair)_{\xi}^{{\rm \gen}}=
\sharp\left\{w\in \aff W^{\lm-\mu}\mid w(\zeta_\pair^\kappa)
=\xi\right\}$ for all $\xi\in \aff\h^*$,
and it is finite if $\kappa\neq 0$.
In particular
${\rm dim}\daff M(\pair)_{\zeta_\pair^\kappa}^{{\rm \gen}}=
\sharp \left( 
\aff W^{\lm-\mu}\cap \aff W[\zeta_\pair^\kappa]\right),$
 where 
$\aff W[\xi]=\{ w\in \aff W
\mid w(\xi)=\xi\}$ for $\xi\in\aff\h^*.$
\end{proposition}
From Proposition~\ref{pr;decom_gen}, it follows that
$\daff M(\pair)$ is an object of $\Rep_\kappa(\daff H_n)$
if $\kappa\in\Z$.
\begin{proposition}\label{pr;piact}
  Let  $(\lm,\mu)\in \In_p$.
Then $\daff M(\pair)\cong \daff M(\varpip\circ(\pair))$.
\end{proposition}
\noindent
{\it Proof.}
Put 
$\lm'=\varpip\circ\lm=
(\lm_p+\kappa-p+1,\lm_1+1,\dots,\lm_{p-1}+1)$
and 
$\mu'=\varpip\circ\mu=
(\mu_p+\kappa-p+1,\mu_1+1,\dots,\mu_{p-1}+1)$.
Put $m=\lm_p-\mu_p$. If $m=0$, then $\daff M(\lm',\mu')\cong 
\daff M(\lm,\mu)$ because
$\zeta_{\lm'\mu'}^\kappa=\zeta_{\pair}^\kappa$ and 
$W_{\lm'-\mu'}=W_{\lm-\mu}.$

Suppose $m\neq 0$.
Set
$v=\pi^{m}{\one_\pair}\in \daff M(\pair)$
(note that $\pi\neq \varpip$).

It can be checked that the weight of $v$
is $\pi^m(\zeta^\kappa_\pair)=\zeta_{\lm',\mu'}^\kappa$, and
that $w v=v$ for $w\in W_{\lm'-\mu'}$.
Hence, there exists a unique $\daff H_n$-homomorphism
 $\psi: \daff M(\lm',\mu')\to
\daff M(\lm,\mu)$
such that $\psi(\one_{\lm',\mu'})=v=\pi^{m}
\one_{\lm,\mu}$.
Similarly, there exists a unique homomorphism
 $\psi': \daff M(\lm,\mu)\to
\daff M(\lm',\mu')$
such that $\psi'
(\one_{\lm,\mu})=\pi^{-m}\one_{\lm',\mu'}$.
Now, it is easy to see that $\psi$ is 
an isomorphism with the inverse $\psi'$.
\qed
\begin{proposition}\label{pr;isomd}
Let $(\lm,\mu)\in \In_p$ and $w\in\aff\Sym_p^\circ$.
Suppose  $w\circ\mu=\mu$. 
Then
 $\daff M(\pair)\cong \daff M(w\circ\lm,\mu).$
\end{proposition}
\noindent
{\it Proof.}
It is enough to prove the statement when 
$w$ is a simple reflection.

Let us prove the statement in the case $w=\sigma_1$ first.
Define
$(\overline{\lm},\overline{\mu})\in \In_2$ and
$(\underline{\lm},\underline{\mu})\in \In_{p-2}$ by
\begin{equation*}
\begin{array}{ll}
\overline\lm=
(\lm_1,\lm_2),\ &\overline\mu=(\mu_1,\mu_2),\\
\underline\lm=
(\lm_3,\dots,\lm_p),\ &\underline\mu=(\mu_3,\dots,\mu_p).
\end{array}
\end{equation*}
Putting $n'=\lm_1-\mu_1+\lm_2-\mu_2$,
 we have
$$\daff M(\pair)\cong \daff H_n(\kappa)
\otimes_{\aff H_{n'}\otimes \aff H_{n-n'}}
\left(\aff M(\overline\lm,\overline\mu)\otimes
\aff M(\underline\lm,\underline\mu)\right).$$

Since $\mu_1-\mu_2+1=0$,
it follows from Lemma~\ref{lem;cpser} that
 $\aff M(\overline\lm,\overline\mu)$ is simple
and $\aff M(\overline\lm,\overline\mu)\cong
\aff M
(\sigma_1\circ\overline\lm,\sigma_1\circ\overline\mu)$.
Hence 
$\daff M(\pair)\cong \daff M(\sigma_1\circ\lm,\mu)$.

Next, 
suppose $\sigma_j\circ\mu=\mu$ $(j\in [0,p-1])$.
Then,
we have 
$\sigma_1\varpip^{1-j}\circ\mu=\varpip^{1-j}\sigma_j\circ\mu=
\varpip^{1-j}\circ\mu$.
Using Proposition~\ref{pr;piact}, we have
\begin{equation*}
\begin{array}{ll}
 \daff M(\pair)&\cong \daff M(\varpip^{1-j}\circ\lm,
\varpip^{1-j}\circ\mu)\\
&\cong\daff M(\sigma_1\varpip^{1-j}\circ\lm,
\varpip^{1-j}\circ\mu)\\
&\cong 
\daff M(\varpip^{j-1}\sigma_1\varpip^{1-j}\circ\lm,\mu)
= \daff M(\sigma_j\circ\lm,\mu).\quad\square
\end{array}
\end{equation*}
\section{Uniqueness of simple quotient}
We give a sufficient condition for 
an induced module $\daff M(\pair)$
to have a unique simple quotient module.

Fix $\kappa\in\Z_{>0}$. Let $\In_p$ and $\zFDaff$ be 
as in \eqref{eq;I_p} and \eqref{eq;I_pk} respectively.
\begin{proposition}
$(\text{\rm{cf. Proposition 2.5.3. in }} \cite{AST})$
\label{pr;onedim}
Let
$\kappa\in \Z_{> 0}$ and $(\lm,\mu)\in \zFDaff$.
Then $\daff M(\lm,\mu)_{\zeta^\kappa_{\lm,\mu}}=
\C\one_{\lm,\mu}$.
\end{proposition}
\noindent
{\it Proof.}
We denote $\zeta^\kappa_{\lm,\mu}$ simply by $\zeta$
till the end of the proof.

It is enough 
to prove the statement in the case $(\pair)\in \FDaff$.

First, suppose $[\mu]_0=\kappa-p+1-(\mu_1-\mu_p)>0$.
Take $u\in \aff W^{\lm-\mu}\cap \aff W[\zeta]$,
where 
$\aff W[\zeta]=\{w\in \aff W\mid w(\zeta)=\zeta\}.$

By Lemma~\ref{lem;representative}, we can write
$u=t_{\eta}\cdot\gm_{\eta}^{-1}
\cdot w= \gm_{\eta}^{-1} \cdot t_{\gm_{\eta}(\eta)} 
\cdot w$
 $(\eta\in P,\,w\in W^{\lm-\mu})$.

Suppose
 $\eta \neq 0$.
Setting $\ech=\sum_{i=1}^n\ech_i$,
we have
$$
\bra\zeta|\ech\ket=
\bra u(\zeta)|\ech \ket
=\bra t_{\gm_{\eta}(\eta)}w(\zeta)|\ech\ket
=\bra \zeta|\ech\ket+\kappa
\bra{\gm_{\eta}(\eta)}|\ech\ket,
$$
and thus $\bra{\gm_{\eta}(\eta)}|\ech\ket=0$.
This implies  that 
$r:=-\bra\gm_{\eta}(\eta)|\ech_1\ket$ is a positive 
integer as
$\gm_{\eta}(\eta)\in P^-$.
Since 
$t_{\gm_{\eta}(\eta)} w(\zeta)=\gm_{\eta}(\zeta)$,
we have
\begin{equation}\begin{array}{rl}\label{eq;contradiction}
0&=\bra t_{\gm_{\eta}(\eta)} w(\zeta)- 
\gm_{\eta}(\zeta)|\ech_1\ket\\
&=\bra w(\zeta)\mid\ech_1\ket
+\bra\gm_{\eta}(\eta)\mid\ech_1\ket\kappa
-\bra\gm_{\eta}(\zeta)\mid\ech_1\ket\\
&=\bra\gm_{\eta}(\eta)|\ech_1\ket\kappa+
\bra\zeta|\ech_{w^{-1}(1)}-\ech_{\gm_{\eta}^{-1}(1)}\ket.
\end{array}
\end{equation}
Put 
$n_0=0$ and $n_i=\sum_{j=1}^{i}(\lm_j-\mu_j)\ (i\in[1,p])
$
as before.
Then 
$w^{-1}(1)=n_{\aa-1}+1$ for some $\aa=[1,p]$ since 
$w\in W^{\lm-\mu}$.
Let $\bb$ be the number such that 
$n_{\bb-1}< \gm_\eta^{-1}(1)\leq n_{\bb}$.
From the definition \eqref{eq;weight_of_one} of $\zeta$, 
it follows that
$\bra\zeta\mid \ech_{w^{-1}(1)}\ket=\mu_\aa-\aa+1$
and $\bra\zeta\mid\ech_{\gm_\eta^{-1}(1)}\ket=
\mu_\bb-\bb+\gm_\eta^{-1}(1)-n_{\bb-1}$.
Now,  (\ref{eq;contradiction})
leads 
\begin{align*}
0&=r\kappa-(\mu_\aa+1-\aa) 
+\mu_\bb-\bb+\gm_\eta^{-1}(1)-n_{\bb-1}\\
&\geq \kappa-(\mu_\aa-\aa)+(\mu_\bb-\bb)\\
&\geq 
\kappa-(\mu_1-1)+(\mu_p-p)
>0.
\end{align*}
This is a contradiction.
Hence  $\eta=0$ and thus $u\in W$.
Therefore
\begin{equation}\label{eq;affW=W}
\aff W^I\cap \aff W[\zeta^\kappa_\pair]=
W^I\cap W[\zeta_\pair].
\end{equation}
This implies
$\daff M(\lm,\mu)_{\zeta^\kappa_\pair}^\gen
=\aff M(\lm,\mu)_{\zeta_\pair}^\gen$
and thus
$$\daff M(\lm,\mu)_{\zeta^\kappa_\pair}=
\aff M(\lm,\mu)_{\zeta_\pair}=\C{\one_\pair}$$
by Proposition~\ref{pr;affine}.

Next, suppose that $[\mu_0]
=0$.
Then there exists $j$ such that $[\mu_j]=\mu_j-\mu_{j+1}+1>0$.
Put 
\begin{equation*}
\lm'=\varpip^{p-j}\circ\lm,\ 
\mu'=\varpip^{p-j}\circ\mu.
\end{equation*}
It is easy to check that $(\lm',\mu')\in\FDaff$
and $[\mu']_0=[\mu]_j>0$.
Moreover we have
$\zeta^\kappa_{\lm',\mu'}=\pi^{n-n_{j}}(\zeta^\kappa_{\lm,\mu})$.
The linear automorphism $v\mapsto \pi^{n-n_{j}}v$ 
on $\daff M(\pair)$ 
gives an isomorphism 
$\daff M(\lm,\mu)_{\zeta^\kappa_{\lm',\mu'}}
{\cong} \daff M(\pair)_{\zeta^\kappa_\pair}.$
On the other hand,
we have an $\daff H_n$-isomorphism
$\daff M(\lm,\mu)\cong\daff M(\lm',\mu')$
by Proposition~\ref{pr;piact}, and thus we have
 $\daff M(\lm,\mu)_{\zeta^\kappa_{\lm',\mu'}}
\cong\daff M(\lm',\mu')_{\zeta^\kappa_{\lm',\mu'}}.$
Therefore,
$\dim \daff M(\lm,\mu)_{\zeta^\kappa_\pair}
=\dim \daff M(\lm',\mu')_{\zeta^\kappa_{\lm',\mu'}}=1$.
\qed
\begin{theorem}\label{th;usqd}
Let $\kappa\in\Z_{>0}$ and $(\pair)\in \zFDaff$.
Then  $\daff  M(\pair)$ has a unique 
simple quotient module, which we denote by $\daff L(\pair)$.
\end{theorem}
\noindent
{\it Proof.} 
Let $N$ be a proper submodule of $\daff M(\pair)$.
By Proposition~\ref{pr;onedim}, we have
$\daff M(\pair)_{\zeta_\pair^\kappa}=\C\one_\pair$.
This implies  $N_{\zeta^\kappa_\pair}=\{0\}$
since ${\one_\pair}$ is a cyclic vector of 
$\daff M(\pair)$.
Hence
$\zeta_\pair^\kappa\notin P(N)$.
Therefore the sum of the all proper submodules of 
$\daff M(\pair)$
is the maximal proper submodule of $\daff M(\pair)$.
\qed

\medskip\noindent
The condition $(\pair)\in\zFDaff$ in  Theorem~\ref{th;usqd} 
can be relaxed by means of Proposition~\ref{pr;isomd}.
\begin{corollary}
\label{co;usqd}
Let  $\kappa\in \Z_{> 0}$. Let $(\lm,\mu)\in \In_p$ and
 $\mu\in\HFDaff$.
Then $\daff M(\lm,\mu)$ has a unique simple quotient module.
\end{corollary}
\section{Classification of simple modules}
Let $\kappa\in\Z_{>0}.$
Let $\Irr(\Rep_\kappa (\daff H_n))$ 
be the set of isomorphism classes of
all simple modules in $\Rep_\kappa (\daff H_n)$.
Through Theorem~\ref{th;usqd}, 
we can construct a correspondence
\begin{equation}\label{eq;Ipktoir}
\tilde\Phi:  \bigsqcup_{p=1}^n
\FDaff
\to \rm{Irr}(\Rep_\kappa(\daff H_n))
\end{equation}
by $(\pair)\mapsto \daff L(\pair)$.
The proofs of the following two theorems are given later.
\begin{theorem}
\label{th;allirrd}
Let $\kappa\in\Z_{>0}$.
Let $\ir$ be a simple module in $\Rep_\kappa (\daff H_n)$.
Then there exists $p\in [1,n]$ and
$(\lm,\mu)\in \FDaff$ such that 
$\ir\cong \daff L(\lm,\mu)$.
In other words, the correspondence $\tilde\Phi$ 
is surjective.
\end{theorem}
\begin{theorem}
  \label{th;isomd}
Let $\kappa\in\Z_{>0}$.
Let $(\pair)\in \FDaff$
and $(\beta,\gamma)\in \In_{q,\kappa}^{*+}$.
Then, the following are equivalent$:$

\smallskip\noindent
${\rm{(a)}}$
$\daff M(\pair)
\cong \daff M(\beta,\gamma)$.

\smallskip\noindent
${\rm{(b)}}$
$\daff L(\pair)
\cong \daff L(\beta,\gamma)$.

\smallskip\noindent
${\rm{(c)}}$
$p=q$ and $(\beta,\gamma)=\varpip^r\circ(\lm,\mu)$
for some $r\in\Z$.
\end{theorem}
By Theorem~\ref{th;isomd} (or Proposition~\ref{pr;piact}), 
the correspondence $\tilde\Phi$ 
factors $\bigsqcup_{p=1}^n \FDaff/\bra \varpip\ket$,
and we get the following.
\begin{corollary}\label{cor;classification}
The correspondence $\tilde\Phi:
\bigsqcup_{p=1}^n \FDaff
\to \rm{Irr}(\Rep_\kappa (\daff H_n))$
above induces a bijection
$$\Phi:
\bigsqcup_{p=1}^n \FDaff/\bra \varpip\ket\to
 \rm{Irr}(\Rep_\kappa (\daff H_n)).$$
\end{corollary}
Recall that $\FDaff\cong \In^*_p/\aff \Sym_p^\circ$ 
(Proposition~\ref{pr;fdpair}) and 
$
\FDaff/\bra \varpip\ket
\cong\In^{*}_p/\aff\Sym_p.$
Hence,
we have a natural one to one correspondence
\begin{equation}\label{eq;Iptoir}
  \bigsqcup_{p=1}^n \In^{*}_p/\aff\Sym_p
\leftrightarrow \rm{Irr}(\Rep_\kappa(\daff H_n)).
\end{equation}
\begin{remark}
We treat the degenerate double affine Hecke algebra 
in this paper, 
but it is easy to modify the arguments 
to obtain the same results for 
the double affine Hecke algebra
provided that 
a certain parameter (often denoted by $q$)
of the algebra is not a root of one. 
In particular, the classification of simple modules over
the double affine Hecke algebra of $\gl_n$ follows.
\end{remark}
\begin{remark}\label{rem;Va}
(i) For (non-degenerate)
double affine Hecke algebras of general type,
a geometric proof of the classification of simple modules
has been given by Vasserot in the preprint~\cite{Va}.

Our parameterization by $(\pair)$  is related to 
Vasserot's parameterization by $\sigma=(\sigma_{a,b})$
in \cite{Va}--\S8 through
$$\sigma_{a,b}=\sharp \{i\in [1,p]\mid
a=\mu_i-i+1,\ b=\lm_i-i\} \ \ (a,b\in\Z,\ a\leq b).$$
(cf. Remark~\ref{rem;quiver}.)

\smallskip\noindent
(ii)
In the preprint~\cite{Ch;mathQA}, 
Cherednik announces a similar classification
for the double affine Hecke algebra of type $A$
by an alternative algebraic approach.
\end{remark}

\begin{remark}
In \cite{AST,Ch;mathQA},
another class of representations have been studied, 
that is, $S(\aff\t)$-semisimple modules.
They form  a subcategory of $\Rep_\kappa(\daff H_n)$,
and the classification 
of simple modules in this category
is given in \cite{Ch;mathQA}.

The method and results developed in our present paper 
are also effective for the study of 
 $S(\aff\t)$-semisimple modules.
We have obtained 
an alternative proof of the 
classification of simple modules
and 
some concrete results
on the structure of simple modules of this class. 
These results will be presented in the forthcoming paper.
\end{remark}
\begin{remark}\label{rem;quiver}
It is known that the set $\bigsqcup_{p=1}^{n}\In_p^{*}/\Sym_p\cong
\Irr\Rep(\aff H_n)$ (Theorem~\ref{th;paraffine})
is naturally indexed by 
isomorphism classes of nilpotent representations of the quiver
of type $A$~\cite{Zel;remarks}.

It can be seen that the the set 
$\bigsqcup_{p=1}^{n}\In_p^{*}/\aff\Sym_p$
above is indexed by isomorphism classes of nilpotent 
representations of the cyclic quiver:

Let $Q_\kappa$ be the quiver of type $A_\kappa^{(1)}$
with the cyclic orientation, i.e. 
the set of vertices is $\Z/\kappa\Z$ and 
the set of morphisms 
consists of the arrows $i\to i+1$ $(i\in\Z/\kappa\Z)$.
Let $S_n$ be the set of isomorphism classes of $n$-dimensional
nilpotent representations of $Q_\kappa$.

Let $Z$ be the set of all pairs $(a,b)$ of integers such that
$a\leq b$ and 
$(a,b)$ is defined up to simultaneous translation by
a multiple of $\kappa$: $(a,b)\sim (a+m\kappa,b+m\kappa),\ m\in\Z$.
It is known that
isomorphism classes of indecomposable finite-dimensional 
representations of $Q_\kappa$ 
are indexed by elements of $Z$,
and
any finite-dimensional nilpotent representation 
of $Q_\kappa$ is decomposed into 
a sum of indecomposable representations
(see \cite{Lu;quivers} for details).
Let $V(a,b)$ denote the indecomposable representation 
corresponding to $(a,b)\in Z$.

Then,
the correspondences
 $\In^*_p\to S_n$ $(p\in[1,n])$ 
defined by
$$(\pair)\mapsto \sum_{i\in [1,p]}V(\mu_i-i+1,\lm_i-i)$$
give rise to
a bijection $\bigsqcup_{p=1}^{n}\In_p^*/\dot\Sym_p \to S_n$.
\end{remark}
\section{Proof of Theorem~\ref{th;allirrd}}
%
In this section, we will give a proof of
Theorem~\ref{th;allirrd}, 
which asserts the surjectivity of
the correspondence $\Phi$ in 
Corollary~\ref{cor;classification}.
For this purpose, we need to introduce some notations.

Fix $\kappa\in\Z_{>0}$.
For an $\daff H_n$-module $N$ and  $(\lm,\mu)\in \In_p$,
set
$$
 N_{[\lm,\mu]}=\{v\in N_{\zeta^\kappa_{\pair}} \mid wv=v 
\hbox{ for } w\in W_{\lm-\mu}\}.
$$
Set 
\begin{align*}
&\BB_p(N)=\{\mu\in \Z^p \mid 
\exists\lm\in \Z^p 
\hbox{ such that }
(\lm,\mu)\in \In^*_p \hbox{ and } N_{[\lm,\mu]}\neq 0\}.
\end{align*}
\begin{example}
If  $p=n$ 
then $\In_n^*=\{(\pair)\in \Z^n\times\Z^n\mid
\lm=(\mu_1+1,\mu_2+1,\dots,\mu_n+1)\}$.
Hence, for $(\pair)\in \In_n^*$, we have
$\zeta_\pair^\kappa=\sum_{i=1}^n 
(\mu_i-i+1)\e_i+\kappa c^*$
and $N_{[\lm,\mu]}$ is nothing but the weight space 
$N_{\zeta_\pair^\kappa}$. 
In particular
$\sharp \BB_n(N)=\sharp P(N)>0$.
\end{example}
For $\mu\in\Z^p$, we put 
\begin{equation*}
[\mu]_0=\kappa-p+1-(\mu_1-\mu_p),\ \ 
[\mu]_i=\mu_i-\mu_{i+1}+1\ \ (i\in [1,p-1])
\end{equation*}
as before.
\begin{lemma}\label{lem;<0}
Let  $\ir$ be a simple module in $\Rep_\kappa(\daff H_n)$ 
and let $p$ be the minimum integer such that 
$\BB_p(K)\neq \emptyset$.
Suppose that  $[\mu]_i<0$ for
 $\mu\in \BB_p(\ir)$ and $i\in[0,p-1]$.
Then $\sigma_i\circ\mu\in \BB_p(\ir)$.
\end{lemma}
\noindent
{\it Proof.}
Let $\mu\in \BB_p(\ir)$.
Let
$\lm$
 be such that
$(\lm,\mu)\in \In^*_p$ and
$\ir_{[\lm,\mu]}\neq 0$.
Put
$n_0=0,\ 
n_j=\sum_{k=1}^{j}(\lm_k-\mu_k)\ \ (j\in[1,p]).
$

First, let us prove the statement when $i\in [1,p-1]$. 
Suppose $[\mu]_i=\mu_i-\mu_{i+1}+1<0$.
Put $\lm_{(i)}=(\lm_i,\lm_{i+1})\in\Z^2$ and 
$\mu_{(i)}=(\mu_i,\mu_{i+1})\in\Z^2$.
Consider the subalgebra
$A:=\aff H_{[n_{i-1}+1,n_{i+1}]}$ of $\daff H_n$,
which we identify with $\aff H_{n_{i+1}-n_{i-1}}$
through Lemma~\ref{lem;Jsub} and Example~\ref{ex;HJ}.

Take $v\in \ir_{[\lm,\mu]}\setminus\{0\}$ and
consider the 
$A$-module $N:=Av$.
Then $N$ is a surjective image of $\aff M(\lm_{(i)},\mu_{(i)})$.

If
$\aff M(\lm_{(i)},\mu_{(i)})$
is irreducible,
then it follows from Lemma~\ref{lem;cpser} that
$N\cong \aff M(\lm_{(i)},\mu_{(i)})\cong 
\aff M(\lm_{(i)}',\mu_{(i)}')$,
where $\lm_{(i)}'=(\lm_{i+1}-1,\lm_i+1)$ and 
$\mu_{(i)}'=(\mu_{i+1}-1,\mu_i+1)$.
Hence there exists 
$a\in A\subset \daff H_n$
such that $av
\in N_{[\lm_{(i)}',\mu_{(i)}']}\setminus\{0\}$.
Clearly, 
$a v\in 
\ir_{[\sigma_i\circ\lm,\sigma_i\circ\mu]}\setminus\{0\}$.
Hence
 $\sigma_i\circ\mu\in \BB_p(\ir)$.

Suppose that $\aff M(\lm_{(i)},\mu_{(i)})$
is reducible.
Then, by Lemma~\ref{lem;cpser},
there exists an exact sequence
$$0\to\aff L(\lm'_{(i)},\mu'_{(i)})\to
\aff M(\lm_{(i)},\mu_{(i)})
\to\aff L(\lm_{(i)},\mu'_{(i)})\to 0.$$
Since $N$ is a surjective image of
 $\aff M(\lm_{(i)},\mu_{(i)})$, 
it is isomorphic to either $\aff M(\lm_{(i)},\mu_{(i)})$ or 
$\aff L(\lm_{(i)},\mu_{(i)}')$.
If $N\cong\aff L(\lm_{(i)},\mu_{(i)}')$, then
we have $\ir_{[\lm,\sigma_i\circ\mu]}\neq 0$.
(Note that $(\lm,\sigma_i\circ\mu)\in \In_p^*$ 
by the assumption of $p$.)

If $N\cong\aff M(\lm_{(i)},\mu_{(i)})$, then
 $N$ contains a
submodule $\aff L(\lm_{(i)}',\mu_{(i)}')$ and thus 
$\ir_{[\sigma_i\circ\lm,\sigma_i\circ\mu]}\neq\ 0$.
In both cases, we have $\sigma_i\circ\mu\in \BB_p(\ir)$.

Next, let us prove the statement for $i=0$.
Suppose $[\mu]_0=\kappa+p-1-\mu_1+\mu_p<0.$

Consider the subalgebra
$A':=\aff H_{[n_{p-1}+1,n+n_1]}$ of $\daff H_n$,
which is identified with $\aff H_{n+n_1-n_{p-1}}$
through Lemma~\ref{lem;Jsub}.
 
Put $\lm_{(0)}=(\lm_p+\kappa-p+1,\lm_1+1)$ and
$\mu_{(0)}=(\mu_p+\kappa-p+1,\mu_1+1)$.
Take $v\in \ir_{[\lm,\mu]}\setminus\{0\}$
 and consider the $A'$-module
$N':=A'v$.
Then $N'$ is a surjective image of the 
$\aff H_{n+n_1-n_{p-1}}$-module
$\aff M(\lm_{(0)},\mu_{(0)})$.
By similar arguments as in the case $i\in [1,p]$,
we have either 
$\ir_{[\sigma_0\circ\lm,\sigma_0\circ\mu]}\neq 0$ 
or $\ir_{[\lm,\sigma_0\circ\mu]}\neq 0$.
Therefore $\sigma_0\circ\mu\in \BB_p(\ir)$.
\qed
\begin{lemma}\label{lem;intersects}
 Let $\kappa\in\Z_{>0}$.
Let $\ir$ be a simple module in $\Rep_\kappa(\daff H_n)$.
If $\BB_p(\ir)\cap \HFDaff\neq \emptyset$,
then there exists $(\pair)\in \FDaff$ such that
$\ir\cong \daff L(\pair)$.
\end{lemma}
\noindent
{\it Proof.}
By the assumption, there exists $(\lm',\mu)\in \In_p^*$ such that 
$\mu\in \HFDaff$ and $\ir_{[\lm',\mu]}\neq 0$.
Take $v\in \ir_{[\lm',\mu]}\setminus\{0\}$.
Since $\ir =\daff H_n v$,  
it is a surjective image of $\daff M(\lm',\mu)$
and thus $K\cong \daff L(\lm',\mu)$ by Theorem~\ref{th;usqd}.
Noting that $\mu\in \HFDaff$,
we can find $w\in \aff\Sym_p$ such that $w\circ\mu=\mu$
and $(w\circ\lm',\mu)\in\FDaff$.
Put $\lm=w\circ\lm'$. 
Now, Proposition~\ref{pr;isomd}
implies $\daff L(\pair)\cong \daff L(\lm',\mu)\cong \ir$.\qed

\medskip\noindent
\underline{{\it Proof of Theorem~\ref{th;allirrd}.}}

\smallskip
Let $K$ be a simple module in $\Rep_\kappa(\daff H_n)$.
Take the smallest integer $p$ such that
$\BB_p(K)\neq\emptyset$.

By Lemma~\ref{lem;intersects}, it is enough to prove that 
$\BB_p(\ir)\cap \HFDaff\neq \emptyset$.
Take $\mu\in \BB_p(K)$, and let $\mu^+$ denote the unique element 
in $\{w\circ\mu\}_{w\in\aff\Sym_p^\circ}\cap \HFDaff$.
Take the shortest $w\in\aff\Sym_p^\circ$ such that
$w\circ\mu=\mu^+$.
Let $w=\sigma_{i_1}\sigma_{i_2}\cdots\sigma_{i_l}$ be 
a reduced expression.
Then by Lemma~\ref{lem;fundshortest},
we have
$[\sigma_{i_{k+1}}\sigma_{i_{k+2}}\cdots 
\sigma_{i_l}\circ\mu]_{i_k}<0$
for $k\in[1,l]$.
Now, Lemma~\ref{lem;<0} implies $\mu^+\in \BB_p(K)\cap\HFDaff$.
\qed
\section{Proof of Theorem~\ref{th;isomd}}
We will give a proof of Theorem~\ref{th;isomd},
which asserts  the injectivity of the correspondence $\Phi$
in Corollary~\ref{cor;classification}.

Fix $\kappa\in\Z_{>0}$.
We start with some preparations.
\begin{lemma}
\label{lem;s_inot}  
Let $(\pair)\in \FDaff$ and
  $s_i\in W_{\lm-\mu}$.
Then $s_i(\zeta_\pair^\kappa)\notin P(\daff M(\pair))$.
\end{lemma}
\noindent{\it Proof.}
We will give a proof only 
for the case  $[\mu]_0>0$.
Other cases can be shown similarly,
using the same argument using 
Proposition~\ref{pr;piact}
as in the proof of Theorem~\ref{th;usqd}.

Assume that there exists  $s_i\in W_{\lm-\mu}$
such that  
$s_i(\zeta_\pair^\kappa)\in 
P(\daff M(\pair))=\aff W^{\lm-\mu}\zeta_\pair^\kappa$
(Proposition~\ref{pr;decom_gen}).
Then $s_i(\zeta_\pair)=x(\zeta_\pair)$
for some $x\in \aff W^{\lm-\mu}$.
Putting $w=s_ix$,
we have $w\in \aff W^{\lm-\mu}\cap\aff W[\zeta_\pair^\kappa]$
because $R(s_ix)=R(x)\sqcup\{x^{-1}(\al_i)\}$ or 
$R(s_ix)=R(x)\setminus\{-x^{-1}(\al_i)\}$.

Recall that
 we have proved
 $\aff W^{\lm-\mu}\cap\aff W[\zeta_{\pair}^\kappa]
= W^{\lm-\mu}\cap W[\zeta_{\pair}]$ 
when $[\mu]_0>0$
in the proof of 
Theorem~\ref{th;usqd} (see \eqref{eq;affW=W}).
Hence $w\in  W^{\lm-\mu}\cap W[\zeta_{\pair}]$.

First, we consider the case where
$$\lm=(m,m+1,\dots,m+p-1),\ \mu=(0,1,\dots,p-1)\text{ with }mp=n.$$
Set $I_j=\{k\in[1,n]\mid \bra \zeta_\pair\mid \ech_k\ket
=\bra \zeta_\pair\mid \ech_j\ket\}$ for $j\in [1,n]$.
Then, in particular, we have
$I_{i+1}=\{k+1\mid k\in I_i\}$.
By induction on $k$, we have
$w(k+1)=w(k)+1$ for all $k\in I_i$.
Taking $k=w^{-1}(i)\in I_i$, we have 
$s_iw(k)=i+1$ and $s_iw(k+1)=i$. 
This implies $x=s_iw\notin \aff W^{\lm-\mu}$.
This is a contradiction.
The same contradiction is deduced 
for general $(\pair)$
through Lemma~\ref{lem;block}.
\qed

\medskip
Let  $\xxi\in P_{\kappa}$.
Then, there exists a unique 
 element $(\xxi^L,\xxi^R)\in\bigsqcup_{p=1}^{n}\In_p^*$
for which the following two conditions hold:
\begin{align}
\label{eq;cLR1}
s_i\in W_{\xxi^R-\xxi^L}
&\Leftrightarrow 
\bra\xxi\mid\al\ch_i\ket=-1
\ \ (i\in[1,n-1]).\\
\label{eq;cLR2}
\zeta_{\xxi^L,\xxi^R}^\kappa&=\xxi.
\end{align}
We denote $h(\xxi)=p$ if $(\xxi^L,\xxi^R)\in \In_p^*$.

In \S 6, we defined
$\zeta_\pair^\kappa\in P_\kappa$ for each
$(\pair)\in\In_p$.
The correspondence $\xxi\mapsto(\xxi^L,\xxi^R)$
is a left inverse of $(\pair)\mapsto \zeta_\pair^\kappa$ in 
the following sense: 
\begin{lemma}\label{lem;op}
If $(\pair)\in\FDaff$ then  
$((\zeta_{\pair}^\kappa)^L,(\zeta_{\pair}^\kappa)^R)=(\lm,\mu)$
and $h(\zeta_\pair^\kappa)=p$.
\end{lemma}
\noindent{\it Proof.}
The statement follows easily from the definition of 
$(\xxi^L,\xxi^R)$.
\qed
\begin{definition}
For a subset $S$ of $P_\kappa$.
Define $\CS(S)$ to be the subset of $S$
consisting of all elements $\xxi\in S$ satisfying 
the following conditions:

\smallskip
\noindent
$\rm{(C1)}\ $ If $\bra \xxi\mid\al\ch_i\ket\ <0$ 
for $i\in[0,n-1]$ then $s_i(\xxi)\notin S$.

\smallskip
\noindent
$\rm{(C2)}\ $ $(\xxi^L,\xxi^R)\in\bigsqcup_{p=1}^n\FDaff.$
\end{definition}
For $(\pair)\in\FDaff$,
put
$\CS_\pair=\CS(P(\daff L(\pair))).$
%
\begin{lemma}\label{lem;zetain}
Let $(\pair)\in \FDaff$. Then $\zeta_\pair^\kappa\in \CS_{\lm,\mu}$.
\end{lemma}
\noindent{\it Proof.}
It is obvious that  $\zeta_\pair^\kappa$ satisfies
(C2) by Lemma~\ref{lem;op}. 
By $(\pair)\in \FDaff$, we have
$$\bra \zeta_\pair^\kappa\mid\al\ch_i\ket\ <0
\Leftrightarrow \bra \zeta_\pair^\kappa\mid\al\ch_i\ket\ =-1
\Leftrightarrow s_i\in W_{\lm-\mu}.$$
Now, it follows from Lemma~\ref{lem;s_inot}
that $\zeta_\pair^\kappa$ satisfies
(C1).\qed

\medskip
We fix $(\pair)\in\FDaff$ for a while.
Let $\xxi\in \CS_\pair$. Put $q=h(\xxi)$ and
put 
\begin{align}\label{eq;n_i}
&n_0=0,\quad 
n_i=\sum_{j=1}^{i}(\lm_j-\mu_j)\ \ (i\in [1,p]),\\
\label{eq;m_i}
&m_0=0,\quad
m_i=\sum_{j=1}^{i}(\xxi^L_j-\xxi^R_j)\ \  (i\in [1,q]).
\end{align}

Noting that $\CS_{\pair}\subseteq P(\daff L(\pair))\subseteq 
P(\daff M(\pair))=\aff W^{\lm-\mu}\zeta_{\pair}^\kappa$,
take
$w\in \aff W^{\lm-\mu}$ such that $\xxi=w(\zeta_{\pair}^\kappa)$.
\begin{lemma}\label{lem;tosimple}
If $\al_i\in R_{\lm-\mu}$ $(i\in[1,n-1])$
then $w(\al_i)=\al_{l}$ for some $l\in[0,n-1]$.
\end{lemma}
\noindent
{\it Proof.}
Let $\al_i\in R_{\lm-\mu}$.
We have
$
w(\e_i)=\e_{j'}+k'\delta
$ and 
$w(\e_{i+1})=\e_j+k\delta$
for some $j,j'\in[1,n]$ and $k,k'\in\Z$.
We have $w(\al_i)=\e_{j'}-\e_j+(k'-k)\delta\in\aff R$ and
\begin{equation}
  \label{eq;jj'}
  -1=\bra \zeta_{\pair}^\kappa\mid \al\ch_i\ket
=\bra \xxi\mid \ech_{j'}-\ech_j\ket +(k'-k)\kappa.
\end{equation}
Since $\al_i\in R_{\lm-\mu}$ and $w\in \aff W^{\lm-\mu}$, we have 
$w(\al_i)\in\aff R^+$.
Note that $w(\al_i)\in\aff \Pi$
if and only if  $j'-j+(k'-k)n=-1$.
We assume that $j'-j+(k'-k)n\neq -1$ and 
will deduce a contradiction.

First, suppose $\brac\xxi{\al\ch_{j-1}}=-1$.
Let $A$ be the subalgebra of $\daff H_n$ generated by
$\{\ech_{j'}+k'c,\,\ech_{j-1}+kc,\,\ech_j+kc,\,
s_{\al'},\, s_{\al_{j-1}} \}$,
where $\al'=\e_{j'}-\e_{j-1}+(k'-k)\delta\in\aff R^+$.
Then, it follows from Lemma~\ref{lem;Jsub}
that 
$A$ is isomorphic 
to the degenerate affine Hecke algebra $\aff H_3$
of $\gl_3$. 
 
Let $v\in \daff L(\pair)_{\xxi}$.
Then we have
\begin{align*}
&{\al}\ch v=0,\ \al\ch_{j-1} v=-v,\  s_{\al_{j-1}}v=v.
\end{align*}
The subspace $Av\subset \daff L(\pair)$ is regarded as an 
$A$-module, and it
a surjective image of the induced module
$\aff M(\udl\lm,\udl\mu)$ over $\aff H_3$
with $\udl\lm=(z+1,z+3),\  \udl\mu=(z,z+1)$,
where $z=\bra\xxi\mid\ech_{j-1}\ket$.
By Lemma~\ref{lem;cpser}, 
$\aff M(\udl\lm,\udl\mu)$ is simple 
and thus $Av\cong\aff M(\udl\lm,\udl\mu)$. 
It follows from 
Lemma~\ref{lem;affine_decom_gen}
that $s_{\al_{j-1}}(\xxi)$ is a
weight of $Av$, and hence
  $s_{\al_{j-1}}(\xxi)=s_{j-1}(\xxi)
\in P(\daff L(\pair))$. 
Combined with
the assumption $\brac\xxi{\al\ch_{j-1}}=-1$,
this contradicts to the condition (C1).

Therefore we must have  $\brac\xxi{\al\ch_{j-1}}\neq -1$.
In this case,
we have $s_{j-1}\notin W_{\xxi^L-\xxi^R}$.
This implies
 $\e_j=
\e_{m_{a-1}+1}$
for some $a\in [1,q]$.
Let $b\in[1,q]$ be the number such that
 ${j'}\in [m_{b-1}+1,m_{b}]$. 
Then, using \eqref{eq;jj'},  we have
\begin{align*}
&\xxi^R_b-b+1 
=\bra \xxi \mid\ech_{m_{b-1}+1}\ket
\leq
 \bra \xxi \mid\ech_{{j'}}\ket\\
&=\bra \xxi \mid\ech_{m_{a-1}+1}\ket-(k'-k)\kappa-1\\
&= \xxi^R_a-a+1-(k'-k)\kappa-1.
\end{align*}
Note that 
 $k'-k\geq 0$ as $w(\al_i)\in \aff R^+$.
If  $b\leq a$, then the inequality
$\xxi^R_b-b+1 <\xxi^R_a-a+1$
contradicts to 
the condition $(\xxi^L,\xxi^R)\in \In_{q,\kappa}^{*+}$.
If $b>a$, then 
we have $k'-k>0$ and 
$(\xxi^R_a-a+1)-(\xxi^R_b-b+1 )\geq 1+(k'-k)\kappa>\kappa$.
This is a contradiction too.
Therefore we have 
$j'-j+(k'-k)n= -1$, and hence $w(\al_i)=\al_{j-1}\in\aff\Pi$.
\qed
\begin{lemma}\label{lem;toend}
$\rm{(i)}$ 
For each $i\in[1,q]$, we have
\begin{align}\label{eq;step2_1}
&w^{-1}(\ech_{m_{i-1}+1}) 
=\ech_{n_{a_i-1}+1}+k_i c,\\
\label{eq;step2_2}
&w^{-1}(\ech_{m_{i}}) 
=\ech_{n_{b_i}}+l_i c
\end{align}
 for some $a_i,b_i\in [1,p]$ and $k_i,l_i\in\Z$.
In particular, we have
\begin{align*}
&\xxi^R_i-i=\mu_{a_i}-a_i+k_i\kappa,\ \
\xxi^L_i-i=\lm_{b_i}-b_i+l_i\kappa\quad (i\in[1,q]).
\end{align*}

\smallskip\noindent
$\rm{(ii)}$
If $q=p$ then $a_i=b_i$ and $k_i=l_i$ 
in \eqref{eq;step2_1}\eqref{eq;step2_2}
for all $i\in[1,p]$.

\smallskip\noindent
$\rm{(iii)}$
The correspondences $[1,q]\to[1,p]$
given by $i\mapsto a_i$ and  $i\mapsto b_i$ 
are injective.
In particular $q\leq p$.
\end{lemma}
\noindent
{\it Proof.}
(i)There exist $r\in[1,n]$ and $k\in\Z$ such that
 $w^{-1}(\ech_{m_{i-1}+1})=\ech_{r}+k c$.
Suppose
$n_{a-1}+1<r\leq n_{a}$ $(a\in[1,p])$.
Then $\al_{r-1}\in R_{\lm-\mu}$ and 
$\brac\xxi{\al\ch_{m_{i-1}}}=\bra\zeta_\pair^\kappa
\mid\al\ch_{r-1}\ket=-1$.
This contradicts to the condition
\eqref{eq;cLR1} and we have \eqref{eq;step2_1}.
Similarly, \eqref{eq;step2_2} follows.

\noindent
(ii) The statement follows easily from Lemma~\ref{lem;tosimple}.

\noindent
(iii) Suppose that 
there exist $i,j\in  [1,q]$ such that
\begin{align*}
&w^{-1}(\ech_{m_{i-1}+1})=\ech_{n_{a-1}+1}+k_i c,\ 
w^{-1}(\ech_{m_{j-1}+1})=\ech_{n_{a-1}+1}+k_j c
\end{align*}
for some $a\in [1,p]$.
Then, 
we have $w^{-1}(\al)=0$ for a root
$\al=\e_{m_{i-1}+1}-\e_{m_{j-1}+1}-(k_i-k_j) \delta$.
This is a contradiction and
thus the correspondence $i\mapsto a_i$ is injective.
\qed

\medskip\noindent
Now, we show the key lemma to the proof of 
Theorem~\ref{th;isomd}.
\begin{lemma}\label{lem;pir}
Let $(\pair)\in\FDaff$. 
Let $\xxi\in \CS_{\lm,\mu}$ and 
suppose that 
$h(\xxi)=
{\rm{max}}\{ h(\zeta)\mid\zeta\in \CS_{\lm,\mu}\}$.
Then 
$h(\xxi)=p$ and
$(\xxi^L,\xxi^R)=\varpip^r\circ(\lm,\mu)$
for some $r\in \Z$.
\end{lemma}
\noindent
{\it Proof.}
We have $h(\xxi)=
{\rm{max}}\{ h(\zeta)\mid\zeta\in \CS_{\lm,\mu}\}=p$ 
by Lemma~\ref{lem;zetain}
 and Lemma~\ref{lem;toend}-(iii).

Take
 $w\in \aff W^{\lm-\mu}$ such that 
$\xxi=w(\zeta_{\pair}^\kappa)$.
By Lemma~\ref{lem;toend}-(i)(ii),
there exist  $a_i\in [1,p]$ and $k_i\in\Z$ such that
\begin{align*}
&w^{-1}(\ech_{m_{i-1}+1})
=\ech_{n_{a_i-1}+1}+k_i c,\ \ 
w^{-1}(\ech_{m_{i}})=\ech_{n_{a_i}}+k_i c
\end{align*}
for each $i\in[1,p]$.

Let  $y$ be the element of $\Sym_p$
such that $y(i)=a_i$ $(i\in[1,p])$, 
and put  $\nu=\sum_{i=1}^p k_i e_i$ ($\nu$ is
an element of the weight lattice of $\gl_p$).
Put 
$x=\tp_\nu y\in\aff\Sym_p$.
Then, we have
$(\xxi^L,\xxi^R)=x\circ(\pair)$.
Since $\FDaff$ is a fundamental domain for the
action of $\aff \Sym^\circ_p$ on $\In_p^*$,
the condition $(\xxi^L,\xxi^R)\in \FDaff$
implies $(\xxi^L,\xxi^R)=\varpip^r\circ(\pair)$
for some $r\in\Z$.
\qed

\medskip\noindent
\underline{{\it Proof of Theorem~\ref{th;isomd}.}}

\noindent
The implication (a)$\Rightarrow$(b) is clear, and
(c)$\Rightarrow$(a) follows from 
Proposition~\ref{pr;piact}.
Let us prove
(b)$\Rightarrow$(c), which 
 completes the proof of Theorem~\ref{th;isomd}.

Suppose $\daff L(\pair)\cong \daff L(\beta,\gamma)$.
Then we have
$P(\daff L(\pair))=P(\daff L(\beta,\gamma))$ and
$\CS_\pair=\CS_{\beta,\gamma}$.
By Lemma~\ref{lem;zetain},
we have $\zeta_{\beta,\gamma}^\kappa\in \CS_{\beta,\gamma}
=\CS_\pair$.
%
By Lemma~\ref{lem;op}, we have 
$((\zeta_{\beta,\gamma}^\kappa)^L,(\zeta_{\beta,\gamma}^\kappa)^R)=
(\beta,\gamma)$ as $(\beta,\gamma)\in\In_{q,\kappa}^{*+}$.
On the other hand, Lemma~\ref{lem;pir} implies
$q=p$ and
$((\zeta_{\beta,\gamma}^\kappa)^L,(\zeta_{\beta,\gamma}^\kappa)^R))
=\varpip^r\circ(\pair)$ for some $r\in\Z$.
Therefore we have $(\beta,\gamma)=\varpip^r\circ(\pair)$.
\qed

\end{document}